\documentclass[graybox,natbib]{svmult}
\usepackage{newtxtext}
\usepackage{helvet}         
\usepackage{courier}        
\usepackage{type1cm}        
\usepackage{makeidx}         
\usepackage{graphicx}        
\usepackage{multicol}        
\usepackage[bottom]{footmisc}
\usepackage{amsfonts}
\usepackage{latexsym}
\usepackage{amsmath}
\usepackage{amssymb}
\usepackage{amsxtra, mathrsfs, upgreek, bm}

\usepackage{bbm}

\makeindex             

\DeclareMathAlphabet{\mathscr}{OT1}{pzc}{m}{it} 
\usepackage{eucal}

\newcommand{\begeq}{\begin{equation}\begin{gathered}}
\newcommand{\eqend}{\end{gathered}\end{equation}}
\newcommand{\begal}{\begin{equation}\begin{aligned}}
\newcommand{\alend}{\end{aligned}\end{equation}}

\begin{document}

\title{  Propagation of chaos for a stochastic particle system modelling epidemics}
\titlerunning{ Propagation of chaos for a stochastic particle system modelling epidemics}
\author{Alessandro Ciallella, Mario Pulvirenti and Sergio Simonella}
\authorrunning{Ciallella, Pulvirenti, Simonella}
\institute{A. Ciallella 
\at  Dipartimento di Ingegneria Civile, Edile -- Architettura e
             Ambientale, and 
             International Research Center M\&MOCS, 
             Universit\'a dell'Aquila,
             via Giovanni Gronchi 18, 67100, L'Aquila, Italy,
            \\
\email{alessandro.ciallella@univaq.it}
\and
M. Pulvirenti
\at Dipartimento di Matematica, Universit\`a di Roma La Sapienza
Piazzale Aldo Moro 5, 00185 Rome, Italy, 
 International Research Center M\&MOCS, Universit\`a dell'Aquila, Italy, and 
 Accademia Nazionale dei Lincei, \\
\email{pulviren@mat.uniroma1.it}
\and
S. Simonella
\at UMPA UMR 5669 CNRS, ENS de Lyon, 46 all\'{e}e d'Italie,
69364 Lyon Cedex 07 -- France,
            \\ 
\email{sergio.simonella@ens-lyon.fr}
}

\maketitle

\abstract{
 We consider a simple stochastic $N$-particle system, already studied by the same authors in \cite{CPS21}, representing  different populations of agents.  Each agent has a label describing his state of health. We show rigorously that, in the limit $N \to \infty$, propagation of chaos  holds, leading to a  set of kinetic equations which are a  spatially inhomogeneous version of the  classical SIR model.
 We improve a similar result obtained in \cite{CPS21} by using here a different coupling technique, which makes the analysis simpler, more natural and transparent.}


\newpage

\section{Introduction\label{cps-intro}}

We consider a stochastic particle system  representing agents moving independently according to a random flight.  
Each agent may be susceptible, infected, or recovered at the initial time.
When an infected particle and a susceptible particle are sufficiently close,  they may interact and become both infected.
Each infected agent can recover independently of the others in a random time of fixed rate. As in the classical SIR model by \cite{KMK},
recovered agents can no longer become susceptible or infected.
However,  at variance with the SIR model, we intend to allow spatially inhomogeneous distributions of the populations. 
 A simple proposal for such a model has been given in \cite{CPS21}, which we reconsider in the present paper. Our purpose is not to provide realistic modelling of spatial patterns (even though this can be important in applications), but rather follow a natural approach, inspired by the kinetic theory of rarefied gases; see for instance \cite {Betal, albi2021kinetic} and references therein. 
There is  therefore no focus on identifying realistic interactions between agents.
In particular, we shall neglect possible individual strategies, and assume a binary interaction.
In spite of its simplicity, the kinetic model still provides a good description of the qualitative behaviour of SIR-like equations, because the essential aspects of the evolution are weakly dependent on the microscopic details (\cite{PS_SIR,CPS20,CPS21}).

In this paper we show that, in the limit $N\to \infty$, the $N$-particle model introduced in \cite{CPS21} reduces to kinetic equations (\eqref{cps:B} below) for the one-particle marginals of the probability measure describing the statistical behaviour of the system. In particular   as a crucial step, we prove {\em propagation of chaos}, namely the asymptotic statistical independence of the agents.
 We use a coupling method, which improves the result obtained in \cite {CPS21} (obtained via the hierarchy method) providing a simpler and more effective proof.
More precisely we introduce a second random process, accounting for the kinetic equations {\em formally} associated with the limit $N \to \infty$, and construct a realization of both processes on the same probability space. This coupling is then used to estimate the distance between the processes and verify that the two models are asymptotically equivalent.
 
The plan of the paper is the following. In Sections \ref{cps:sec:model} and \ref{cps:sec:K_L} we present the model and the limiting kinetic equation, respectively. 
In Section \ref{cps:sec:P_A} we discuss the convergence of the particle model in the kinetic limit  following the coupling argument.
Finally, Section \ref{cps:sec:conclusion} is devoted to concluding remarks.

\section{Model}\label{cps:sec:model}

We now define the system we are going to study. We refer also to \cite{CPS21} for the description of the model and the derivation of the formal kinetic limit. 
Consider $N$ particles, representing the agents of the system, moving 
a square in the plane $\mathbb{R}^2$ with periodic boundary conditions, i.e., the torus $\Lambda=[0,D]\times[0,D]$, where $D>0$ is a fixed parameter.
The particles are assumed to move with velocities of modulus $1$ so that the velocities belong to the unit circle $\mathcal{S}^1$.
Each particle has a label $a_i\in \{S,I,R\}=:L$, $i=1,2,\ldots,N$, representing the class of the agent (susceptible, infected, recovered).
We introduce the notations $Z_N=(z_1,z_2,\ldots,z_N)$, where $z_i=(x_i,v_i)\in \Lambda\times \mathcal{S}^1$, and $A_N=(a_1,a_2,\ldots,a_N)$.
At any time, the state of the system is described by the  $\left(4N\right)$--tuple of coordinates  $(Z_N;A_N)$ belonging to the phase space $ \left(\Gamma\times L\right)^N$, where $\Gamma:=\Lambda\times \mathcal{S}^1$.

The system evolves in time according to the following generator:
\begin{equation}
\label{cps:gen1}
 {\mathcal L}= {\mathcal L}_0+{\mathcal L}_1+ {\mathcal L}_d +{\mathcal L}^N_{int},
\end{equation}
where we have generators for which we omit the dependence on $N$ since they are acting independently on single particles ($\mathcal L_0$, $\mathcal L_1$ and $\mathcal L_d$),  and a binary interaction term between particles ${\mathcal L}^N_{int}$. 

More precisely, ${\mathcal L}_0+{\mathcal L}_1$ is the generator of $N$ independent copies of a random walk:
\begin{equation} \label{cps:eq:L0}
\mathcal{L}_0= \sum_{i=1}^N v_i \cdot \nabla_{x_i}
\end{equation}
 is the generator of free motion (particles are moving of linear motion with velocities $v_i$, $i=1,\ldots,N$), and jumps happen as described by
\begin{equation*}
\mathcal{L}_1 \Phi(Z_N) = \sum_{i=1}^N \frac{1}{2\pi} \int_{\mathcal{S}^1} dw [\Phi(z_1,\ldots,x_i,w,\ldots,z_N)-\Phi(Z_N)].
\end{equation*}
In particular, when we select the test function $\Phi$ of the form of a function of a single particle state only, for instance of the first one $\Phi(Z_N)=\phi(z_1)$, we get
\begin{equation} \label{cps:eq:L1}
\mathcal{L}_1 \Phi(Z_N) =  \frac{1}{2\pi} \int_{\mathcal{S}^1} dw [\phi(x_1,w)-\phi (x_1,v_1)].
\end{equation}
Note that labels are not involved in the random flight process, so here we have left out the dependency on them to shorten the notations.

The other contributions in the generator are instead acting on the labels of the particles only.
The term $\mathcal{L}_d$ describes the decay of infected agents $I$ into recovered $R$. It takes the form
\begin{equation} \label{cps:eq:Ld}
\mathcal{L}_d \Phi(Z_N;A_N)= \gamma \sum_{i=1}^N [\Phi(Z_N;a_1,\ldots,\tilde{a}_i,\ldots,a_N)-\Phi(Z_N;A_N)],
\end{equation}
$\gamma>0$ being a constant parameter representing the rate of decay, and where the transition of an $a_i$ label into a $\tilde{a}_i$ is defined by
\begin{equation}
\label{cps:eq:decay_def}
\tilde a_i=R \,\,\,\,\,\, \text {if} \,\,\,\,\,a_i=I; \qquad \tilde a_i=a_i  \,\,\,\,\, \text{otherwise}.
\end{equation}

Finally, the binary interaction term describing the infection process has generator
\begin{equation*}
\begin{split}
\mathcal{L}_{int}^N \Phi(Z_N;A_N)= \frac{\lambda}{N} \sum_{\substack {i,j =1 \\ j> i}}^N \big[ &
\Phi(Z_N;a_1,\ldots ,a_i',\ldots ,a_j',\ldots ,a_N)
\\ &  - \Phi(Z_N;a_1,\ldots ,a_N)\big].
\end{split}
\end{equation*}
Here, $\lambda>0$ is the constant representing the rate of the process and for the involved particles the transition from $(a_i,a_j)$ into $(a_i',a_j')$ is defined by 
\begin{equation}
\label{cps:a'}
\begin{cases}
\text{if} \,\,\, \chi_{i,j}=1   \,\,\, \text{and} \,\,\, a_i=I, a_j =S \,\,\, \text {or} \,\,\, a_j=I, a_i =S \,\,\, \text {then} \,\,\,a'_i=a'_j=I \,  ;   \\
  \text {otherwise} \,\,\,a_i'=a_i, \, a_j' =a_j .
\end{cases}
\end{equation}
The characteristic function $\chi_{i,j}$ is introduced to allow two particles to interact only when they are sufficiently close: fixing $R_0>0$ we define
\begin{equation*}
\chi_{i,j}:=\mathbbm{1}_{\{ x_i,x_j \ |\ |x_i-x_j | <R_0 \}}\;.\
\end{equation*}

Thus, the evolution due to the defined generator \eqref{cps:gen1} is describing the following behaviour. Each agent is moving in the space $\Lambda$ performing a random flight, where the velocity jumps happen with rate $1$.
Each particle has a label describing its state: susceptible (S), infected (I), or recovered (R).
According to a Poisson process with overall rate 
$\frac {(N-1)\lambda}2$, a pair of agents is selected uniformly. Whenever the chosen agents are at a distance smaller than $R_0$, the infection process takes place only for pairs made by an infected and a susceptible agent, which are changed in a couple of infected agents.
Finally, each infected agent becomes recovered according to a different Poisson process with rate $\gamma$.

A statistical description is in order for dealing with such a system with a large number of agents $N$.
The initial configuration of the system at time zero is given by the probability density $W_0^N$, \emph{symmetric} in the exchange of particles, such that
\begin{equation*}
W_0^N:(\Gamma\times L)^N \to \mathbb{R}^+,\qquad \qquad \sum_{A_N} \int_{\Gamma^N} d Z_N W_0^N (Z_N;A_N)=1.
\end{equation*}
The time evolved measure $ W_t^N (Z_N;A_N), t>0$ is given by
\begin{equation*}
\begin{split}
&\sum_{A_N} \int dZ_N W_t^N (Z_N;A_N) \Phi (Z_N;A_N)= \\
&\sum_{A_N} \int dZ_N W_0^N (Z_N;A_N)  \mathbb E  [\Phi (Z_N(t);A_N(t) )]\;,
\end{split}
\end{equation*}
where $\Phi$ is a test function, $(Z_N;A_N)  \to (Z_N(t);A_N(t))$ is the process and $\mathbb E=\mathbb E_{(Z_N,A_N)} $ is the expectation conditioned to the initial value $(Z_N;A_N)$.
Integrating with respect to the last $N-j$ particle positions, velocities and labels, $j=1,\ldots,N-1$, we build the $j$-particle marginal $f_j^N$, that gives the probability density of finding $j$ particles with labels $A_j$ in the configuration $Z_j$:
\begin{equation*}
f^N_j(Z_j;A_j;t)= \sum_{\bar A_{N-j} \in L^{N-j}} \int d\bar Z_{N-j} W_t^N ( Z_j,\bar Z_{N-j};A_j,\bar A_{N-j})\;.
\end{equation*}

We consider particles distributed independently at time zero. 
Given the one-particle density distribution $f_0$, with normalization
\begin{equation} \label{cps:eq:norm0}
\sum_{a\in L} \int dz f_0(z;a)=1\;,
\end{equation}
the initial state is
\begin{equation} \label{cps:eq:fullindtz}
W_0^N (Z_N;A_N)=\prod_{i=1}^N f_0(x_i,v_i; a_i):= (f_0 )^{\otimes N} (Z_N;A_N).
\end{equation}

It is crucial to keep in mind that, even when at time zero the particles are independent, they do not remain independent at positive times.
In fact, the dynamics creates correlations between particles so that the measure is no longer factorized at positive times.
However, we will prove that the so--called \lq propagation of chaos\rq holds. This means that this independence is recovered asymptotically in the limit $N\to \infty$.

\section{Kinetic limit}\label{cps:sec:K_L}

We aim at proving that the system we introduced in the previous section is asymptotically equivalent to the system one  formally obtains in the limit $N\to \infty$ assuming the propagation of chaos.
We refer to \cite{CPS21} for such   formal derivation, and for a proof of the rigorous result through the hierarchy of equations satisfied by the marginals $f_j^N(Z_j;A_j;t)$.
 In the following section, we shall proceed in a simpler and more natural way.

According to the formal limit, the triple of single-particle densities $$\left(f(z;S;t), f(z;I;t), f(z;R;t)\right)$$ satisfies the following system of  kinetic equations ($z=(x,v)$):
\begin{equation}
\label{cps:B}
\begin{cases}
  \left(\partial_t + v \cdot \nabla_x\right) f(z;S) =&  {\mathcal L}_1 f (z;S) - \lambda  f(z;S) \int f(z_1;I) \chi (|x-x_1| <R_0) dz_1 \\
    \left(\partial_t + v \cdot \nabla_x\right) f(z;I)  =& {\mathcal L}_1f(z;I)  -\gamma f(z;I)
    \\ &\qquad+\,\,\lambda  f(z;S)  \int f(z_1 ;I) \chi (|x-x_1| <R_0) dz_1   \\
     \left(\partial_t + v \cdot \nabla_x\right) f (z;R) =& {\mathcal L}_1 f(z;R ) + \gamma f(z;I) 
\end{cases}\, .
\end{equation}
Note that the sum
$$
f(z,t) := \sum_{a \in L} f(z;a;t)
$$
satisfies the simple random flight equation
$$
\left(\partial_t + v \cdot \nabla_x\right) f(z,t) =  {\mathcal L}_1 f (z,t)\;.
$$
Moreover, the system of kinetic equations \eqref{cps:B} provides
a more detailed description of the classical SIR model when dealing with  spatially inhomogeneous data
(see \cite{CPS21}).

We now introduce as second model the process associated with the kinetic equations \eqref{cps:B} that we will eventually couple with the  one introduced in Section \ref{cps:sec:model}.
 The single-particle generator is
\begin{equation*}
\label{cps:gen2}
\tilde {\cal L}= \tilde {\cal L}_0+\tilde {\cal L}_1+\tilde {\cal L}_d +{\cal \tilde L}_{i}
\end{equation*}
 where $\tilde{\cal L}_0$, $\tilde{\cal L}_1$ and $\tilde{\cal L}_d$ are the operators defined in Section \ref{cps:sec:model}, namely ${\cal L}_0$, ${\cal L}_1$ and ${\cal L}_d$ respectively, in the case $N=1$; while the infection generator is given by the non--linear ($f$-dependent) term
\begin{equation} \label{eq:genBp}
{\cal \tilde L}_{i} \phi(z;b)= \lambda \mathcal{N}_f(z)\{  \phi (z;b') - \phi (z;b) \},
\end{equation}
where
$\mathcal{N}_f $  (depending on time through $f$) is defined as
\begin{equation*}
\mathcal{N}_f (z) =(f(I;t) * \chi_{R_0}) (z) := \int f(z_1;I;t) \chi (|x-x_1| <R_0) dz_1\;,
\end{equation*}
$\chi(A)$ is the indicator function of $A$, and $(z;b)\in \Gamma\times L$. 
We are keeping the same notation associated with the prime sign  as in Section \ref{cps:sec:model}, i.e., $b'=I$ if $b=S$, while in the other cases $b'=b$.

We now consider the $N$--particle process defined by $N$ independent copies of the one--particle nonlinear process introduced above.
The  $N$--particle generator is
\begin{equation*}\label{cps:gen_tilde}
\tilde{\cal L}^N= {\cal L}_0+{\cal L}_1+ {\cal L}_d +{\cal \tilde L}_{i}^N 
\end{equation*} 
 where ${\cal L}_0$, ${\cal L}_1$, ${\cal L}_d$ are given by \eqref{cps:eq:L0}-\eqref{cps:eq:Ld} and,  for $(Z_N;B_N)\in (\Gamma\times L)^N$,
\begin{equation*}
{\cal \tilde L}_{i}^N \Phi(Z_N;B_N)=\lambda \sum_{i=1}^N \mathcal{N}_f(z_i) \left\{  \Phi(Z_N;B_N^i)-\Phi(Z_N;B_N)  \right\}
\end{equation*}
with 
\begin{equation}
\label{cps:eq:apex_i'}
\begin{cases}B_N^i\!=\!(b_1,\ldots,b_{i-1},I,b_{i+1},\ldots,b_{N}) \,\,\,
\text{if} \,\,\, B_N\!=\!(b_1,\ldots,b_{i-1},S,b_{i+1},\ldots,b_{N})   \\
B_N^i=B_N \,\,\, \text {otherwise} \\
\end{cases}.
\end{equation}

\section{Particle approximation}\label{cps:sec:P_A}

We are now ready to prove the  validity of the kinetic equations \eqref{cps:B} for the model introduced in Section \ref{cps:sec:model}, in the limit $N \to \infty$. We are assuming \eqref{cps:eq:norm0}  and \eqref{cps:eq:fullindtz} at time zero, i.e., the initial datum is factorized.
On the other hand the dynamics of labels creates correlations as
$$
f^N_2 ( z_1,z_2; a_1,a_2; t) \neq f^N_1 ( z_1;a_1; t) f^N_1 (z_2;a_2; t)\;.
$$
To show that these correlations are vanishing as $N \to \infty$,
we now introduce a coupling between the process described by the generator \eqref{cps:gen1} and the process associated to the kinetic equations \eqref{cps:B}. 

We start by considering a different equivalent form for the generator
\begin{equation*}
\begin{split}
& {\mathcal L}^N_{int}\Phi  (Z_N;A_N ) 
\\
& =
 \frac{\lambda}{2N} \sum_{i=1}^N\sum_{\substack {j =1 \\ j\neq i}}^N
\delta_{a_j,I} \{	\Phi(Z_N;a_1,\ldots ,a_i',\ldots ,a_j',\ldots ,a_N)
 - \Phi(Z_N;A_N)  \}
\\ &=  \frac{\lambda}{N} \sum_{i=1}^N\sum_{\substack {j =1 \\ j\neq i}}^N
\delta_{a_j,I} { \chi_{i,j}} \{	\Phi(Z_N;A_N^i) - \Phi(Z_N;A_N)   	\}
\\&=
 \lambda
\sum_{i=1}^N
\mathcal{J}_{emp}^i \{\Phi(Z_N;A_N^i) - \Phi(Z_N;A_N)\} \, ,
\end{split}
\end{equation*}
where we use the prime superscript with the meaning of \eqref{cps:a'}, i.e
\begin{equation}
\nonumber
\begin{cases}a'_i=I \,\,\,
\text{if} \,\,\,  a_j=I, a_i =S  \,\,\, \text{and} \,\,\, \chi_{i,j}=1, \,\,\,  \\
 a_i'=a_i, \, a_j =a'_j \,\,\, \text {otherwise} \\
\end{cases}\, ,
\end{equation}
and the $i$ superscript  (as in \eqref{cps:eq:apex_i'}) when we prefer to make explicit the role of both $\delta_{a_j,I}$ and $  \chi_{i,j}$, i.e.
\begin{equation}
\label{cps:eq:apex_i}
\begin{cases}A_N^i\!=\!(a_1,\ldots,a_{i-1},I,a_{i+1},\ldots,a_{N}) \,\,\,
\text{if} \,\,\, A_N\!=\!(a_1,\ldots,a_{i-1},S,a_{i+1},\ldots,a_{N})   \\
A_N^i=A_N \,\,\, \text {otherwise} \\
\end{cases}
\end{equation}
having introduced 
\begin{equation}
\nonumber
\mathcal{J}^i_{emp} =\mathcal{J}^i_{emp}(Z_N;A_N)=\frac{1}{N}\sum_{\substack {j =1 \\ j\neq i}}^N \delta_{a_j,I} { \chi_{i,j}}
.
\end{equation}

In words, we are interpreting the behaviour in the following equivalent way. 
The $N$ agents 
are still moving in $\Lambda$ via a random flight, with velocity jumps in $\mathcal{S}^1$ taking place with rate $1$ per agent and each infected agent becomes recovered according to a Poisson process of rate $\gamma$. 
In the infection process we choose with rate $N$ an agent  $ i$ uniformly and then, with a rate $\lambda$, we evaluate the interaction with an uniformly picked agent.
 So that, if this latter is an infected agent and the $i$-th is a susceptible one and if their distance is smaller than $R_0$, then $a_i$ becomes $I$, otherwise nothing happens.

We want to compare the behaviour of this process with that of the one--particle non--linear process associated to the kinetic
equation \eqref{cps:B}, described by the generator  \eqref{eq:genBp}.
The idea is now to couple the two processes. We remark that the positions and velocities are the same, so that the difference is only in the distributions of the labels.
Let us call $\mathcal{E}^i$ the following difference 
\begin{equation*}
\mathcal{E}^i:=\mathcal{N}_f(z_i)-\mathcal{J}^i_{emp} (Z_N;B_N).
\end{equation*}
We expect that, by the law of large numbers, $\mathcal{E}^i$ gives a vanishing contribution as $N\to \infty$ and hence, for sufficiently large $N$, generates a small perturbation in the evolution of distributions of labels.

We consider as coupling the process  $t\to (Z_N(t); A_N(t),B_N(t))$ with generator
\begin{equation*}
Q^N=\mathcal{L}_0+ \mathcal{L}_1+ Q^N_d + Q^N_i,
\end{equation*}
where the decay term is defined as
\begin{equation}
\label{cps:Q^N_d_def}
Q^N_d \Phi(Z_N;A_N,B_N)= \gamma \sum_{i=1}^N [\Phi (Z_N;\tilde{A}^i_N,\tilde{B}^i_N)-\Phi(Z_N;A_N,B_N)],
\end{equation}
denoting by $\tilde{A}^i_N=(a_1,\ldots,\tilde{a}_i,\ldots,a_N)$ the transition on the $i$--th agent as defined in \eqref{cps:eq:decay_def}, and analogous definition for $\tilde B_N^i$.
The interaction term is 
\begin{equation}
\label{cps:Q_N_i_def}
\begin{split}
Q^N_i \Phi(Z_N;A_N,B_N)=\,\, &  \lambda 	\sum_{i=1}^N \mathcal{J}^i \{\Phi(Z_N;A^i_N,B_N^i)-\Phi(Z_N,A_N,B_N)\}
\\  +& \lambda \sum_{i=1}^N \mathcal{J}^i_1 \{\Phi(Z_N;A^i_N,B_N)-\Phi(Z_N;A_N,B_N)\}
\\  +&  \lambda  \sum_{i=1}^N \mathcal{J}^i_2 \{\Phi(Z_N;A_N,B_N^i)-\Phi(Z_N;A_N,B_N)\}
\\  +& \lambda  \sum_{i=1}^N \mathcal{E}^i \{\Phi(Z_N;A_N,B_N^i)-\Phi(Z_N;A_N,B_N)\},
\end{split}
\end{equation}
where $A^i_N$, $B_N^i$ follows the transition rule \eqref{cps:eq:apex_i}, and the label jump  associated to  $\mathcal{J}^i_1$  only affects the transition for $A_N$, while  the terms corresponding to $\mathcal{J}^i_2$ and $\mathcal{E}^i$  only affect the transition for $B_N$. Instead, the jumps described by the first
term in \eqref{cps:Q_N_i_def} and in \eqref{cps:Q^N_d_def} are simultaneous for particles labelled by  both $A_N$ and $B_N$  (which is optimal for what concerns the propagation of chaos).
 We set
\begin{equation}
\nonumber
\begin{split}
&\mathcal{J}^i =\frac{1}{N}\sum_{\substack {j =1 \\ j\neq i}}^N \delta_{a_j,I} \delta_{b_j,I} { \chi_{i,j}},\,\,
\\&\mathcal{J}^i_1 =\frac{1}{N}\sum_{\substack {j =1 \\ j\neq i}}^N \delta_{a_j,I}(1- \delta_{b_j,I}) { \chi_{i,j}},\,\,
\\&\mathcal{J}^i_2 =\frac{1}{N}\sum_{\substack {j =1 \\ j\neq i}}^N (1-\delta_{a_j,I}) \delta_{b_j,I} { \chi_{i,j}}
\end{split}
\end{equation}
(being $\delta_{a_j,I}$ the Dirac delta taking value $1$ if and only if $a_j=I$ and $0$ otherwise, and analogously for the other deltas).
Note that thanks to this choice, one has that
\begin{equation}
\nonumber
\mathcal{J}^i+\mathcal{J}^i_1=\mathcal{J}_{emp}^i(Z_N;A_N),\qquad \mathcal{J}^i+\mathcal{J}^i_2+\mathcal{E}^i=\mathcal{J}^i_{emp}(Z_N;B_N)+\mathcal{E}^i= \mathcal{N}_f(z_i)
\end{equation}
so that $Q^N$ is actually generating a coupling of the two previously described processes, i.e., we recover as the two marginals the two processes we are considering, obtained by integrating over $B_N$ and $A_N$, respectively.

Let $R^N(t)$ be the law at time $t$ for the coupled process, and
we consider as the initial distribution at time $0$ the factorized distribution with the same distribution of labels for
both $a$ and $b$--labelled particles, i.e. $R^N(0)=f_0^{\otimes N}(Z_N;A_N)\delta_{A_N,B_N}$.
We define $D_N(t)$  as the average fraction of particles having different labels $a_i$, $b_i$, i.e.
we choose the test function $\Phi(Z_N;A_N,B_N)=\frac{1}{N} \sum_{i=1}^N d(a_i,b_i)$, where $d(a,b)=1-\delta_{a,b}$:
thanks to the symmetry
\begin{equation*}
D_N(t)=\int d R^N(t) \frac{1}{N} \sum_{i=1}^N d(a_i,b_i)=
\int d R^N(t)  d(a_1,b_1).
\end{equation*}
We notice that $D_N(0)=0$ and that $D_N(t)$ is positive. 
Our aim is to prove that $D_N(t)$ can be estimated by an arbitrarily small quantity  for finite time, provided that $N$ is sufficiently large.

\begin{theorem} \label{cps:thm:main}
 For any $t>0$, the fraction of particles with different labels satisfies 
\begin{equation*}
D_N(t) \leq  \frac {t\lambda} {N} e^{2\lambda t} \;.
\end{equation*}
for sufficiently large $N$.
\end{theorem}

{ \it Remark.  } The above estimate says that, for any positive time, for the time evolved joint measure associated to the $N$--particle systems of the process whose generator is \eqref{cps:gen1}, and of the one whose generator is  \eqref{eq:genBp}, the fraction of particles with different labels is 
arbitrarily small, provided that $N$ is sufficiently large.

\bigskip
{\em Proof.}
We evaluate
\begin{equation}
\begin{split}
\frac{d}{dt} \int d R^N(t) d(a_1,b_1)=\int d R^N(t) (\mathcal{L}_0+ \mathcal{L}_1+ Q^N_d +  Q^N_1) d(a_1,b_1)\;,
\end{split}
\end{equation}
where the contribution due to the terms $\mathcal{L}_0$ and $\mathcal{L}_1$ is vanishing since $d(a_1,b_1)$ does not depend on $x$ and $v$.

The contribution due to the decay term \eqref{cps:Q^N_d_def} is
\begin{equation*}
\gamma \int dR^N(t) (d(\tilde{a}_1,\tilde{b}_1) - d(a_1,b_1))= -\gamma \int dR^N(t) (\delta_{\tilde{a}_1,\tilde{b}_1} -  \delta_{a_1,b_1}) .
\end{equation*}
It is easy to verify that this contribution is non--positive. Indeed, in the case where $a_1$ and $b_1$ are equal, also $\tilde a_1$ and $\tilde b_1$ remain equal, and the contribution is vanishing.
In the case when $a_1$ and $b_1$ are different, after the transformation $\tilde a_1$ and $\tilde b_1$ can either remain different, or (if at the beginning there were an $I$ and an $R$--labelled agent) both become $R$.
This last case is the only one in which the contribution is not vanishing, but negative.

The same argument applies to the first term in the right hand side of the interaction \eqref{cps:Q_N_i_def}. This term can be written as
\begin{equation*}
-  \lambda  \int dR^N(t) \mathcal{J}^1 (\delta_{{a}'_1,{b}'_1} -  \delta_{a_1,b_1}) ,
\end{equation*}
and, in the same way as previously, $a_1=b_1$ implies that after the transformation $a_1'=b_1'$, so that the only non--zero contribution is negative, and it happens when an $S$ and an $I$--labelled particle are transformed into two $I$--particles.

Now, let us evaluate the second term in the right hand side of the interaction \eqref{cps:Q_N_i_def}.
It can be written as
\begin{equation*}
\begin{split}
\mathrm{(II)}=&
  \lambda \int dR^N(t) \mathcal{J}^1_1 (Z_N; A_N,B_N) ( d(a_1',b_1)- d(a_1,b_1) ) 
\\=&
-  \lambda  \int dR^N(t) \mathcal{J}^1_1 (Z_N; A_N,B_N) ( \delta_{a_1',b_1} - \delta_{a_1,b_1} ).
\end{split}
\end{equation*}
Here, the difference of the deltas can take the values -1,0, or 1. Then
\begin{equation}
\begin{split}
|\mathrm{(II)}| &\leq   \lambda  \int dR^N(t) \frac{1}{N} \sum_{ {j =2}}^N \delta_{a_j,I}(1- \delta_{b_j,I}) { \chi_{1,j}}
\\ & \leq  \lambda \int dR^N(t) \frac{1}{N} \sum_{j=2}^N d(a_j,b_j)\leq {\color{cyan} \lambda} D_N(t).
\end{split}
\end{equation}
The third term in \eqref{cps:Q_N_i_def} can be estimated in exactly the same way. 

Finally, 
\begin{equation*}
\begin{split}
\int dR^N(t) \mathcal{E}^1 =&\int dR^N(t) \mathcal{N}_f(z_1) - \int dR^N(t) \mathcal{J}_{emp}^1(Z_N;B_N)
\\
=&\sum_a \int dz \,dz_1 f(z;a;t) \chi(|x-x_1| < R_0) f(z_1;I;t) 
\\&- \int d R^N(t) \frac{1}{N} \sum_{j=2}^N \delta_{b_j,I} \chi(|x_j-x_1|< R_0)\;.
\end{split}
\end{equation*}
The last term, thanks to the symmetry, can be written as 
\begin{equation*}
-\frac{N-1}{N} \sum_a \int dz\, dz_1 f(z;a;t) \chi(|x-x_1| < R_0) f(z_1;I;t),
\end{equation*}
so that
\begin{equation*}
\left|\int dR^N(t) \mathcal{E}^1\right| \leq \frac 1N \sum_a \int dz\, dz_1 f(z;a;t) \chi(|x-x_1|<R_0) f(z_1;I;t) \leq \frac 1N\;.
\end{equation*}

In conclusion, recalling that $D_N(0)=0$, we find that
\begin{equation}
\frac{d }{dt}D_N(t) \leq   2 \lambda D_N(t) + \frac{\lambda}{N}
\end{equation}
and the Theorem is proven by Gr\"{o}nwall's Lemma.
\qed


\bigskip

We conclude this section with some comments.

Theorem 1 implies the convergence of the one-particle marginal $f^N_1(z;a;t) $ to the solution $f( z;a;t) $ of the kinetic equations \eqref {cps:B}. To see this  one can use the Wasserstein distance $\overline{\cal W}$, which we recall. Given two measures $\mu$ and $\nu$ on a metric space $X$ with the discrete metric  $\bar d$,
$$
\overline{\cal W} (\mu, \nu) := \inf_{R \in {\cal C}(\mu,\nu)} \int dR(z,z') \bar d(z,z')
$$
where ${\cal C}(\mu,\nu)$ is the set of the joint representations of $\mu$ and $\nu$, i.e.  measures on $X \otimes X$
with marginals $\mu$ and $\nu$ respectively.  This metric is equivalent to the Total Variation and hence to the $L^1$--norm, whenever such measures have suitable densities  (\cite{villani2003topics}).

Arguing now on the phase space $\Gamma\times L$ endowed with the discrete metric
$$
\bar d(z_1,z_2) +d(a_1,a_2),
$$
we construct a joint representation of $df^N_1(z;a;t) $ and $df( z;a;t) $, denoted by  $dR^N(t)$,  along the coupled process.
 Since the trajectories  (on $\Gamma$) of the two processes are the same, the configurational contribution vanishes and
$$
\overline{ \cal W} \left(df^N_1(t), df(t) \right) \leq  
\frac {t\lambda} {N} e^{2\lambda t}\;.
$$
This also controls $\| f^N_1(t)- f(t) \|_{L^1} $  and, along similar lines, on can also show that
$$
\| f^N_j(t)- f(t)^{\otimes j} \|_{L^1} \leq Cj\frac {t\lambda} {N} e^{2\lambda t}
$$
for some geometrical positive constant $C$.

The obtained estimates slightly improve the quantitative results obtained via the hierarchy
as sketched  in \cite{CPS21}.

\section{Concluding remarks\label{cps:sec:conclusion}}

In this paper we performed a mathematical analysis, based on a coupling method, of a toy model for epidemic spread in a system of $N$ individual agents as proposed in \cite{CPS21}.

 The study of epidemiological models has been obviously very fertile over recent times (see  \cite{PS_SIR,FGPR20,CPS20,bertaglia2021spatial,boscheri2021modeling,CPS21,loy2021viral,VHP21} for a very non-exhaustive list of recent mathematical works oriented toward kinetic theory). Even in this context, the use of coupling methods is not new. We can mention classical work as \cite{whittle55} on the comparison of birth-death processes with the SIR model, 
and several other contributions 
over more recent times
as for instance \cite{ball95,hagg98}.  Of course
coupling arguments have been a powerful tool for studying the behaviour of stochastic processes for several years now. 
In particular such arguments have been used to approximate kinetic equations of Boltzmann type (and justify related numerical schemes) with particle systems in the spirit of Kac; see e.g.\,\cite{graham1997stochastic,CF16}.
 
 Multi-agent systems are important in a wide range of fields. As an example in the theory of epidemic spread, the most common models divide the population in several classes, going from the classical susceptible-infected-recovered quoted above, to recent generalizations with several additional species (as the SUIHTER model with susceptible uninfected individuals, undetected infected, isolated, hospitalized, threatened, extinct, and recovered, proposed in \cite{multiSIR21}). The models aim at describing the evolution of average fractions for the different species, as in the predator-prey system of Lotka-Volterra type (see \cite{LV21} for a recent review). We would like to mention here, as another example, the case of mathematical models describing the behaviour of bones, which include cellular automata models (for instance \cite{cBCA,vsBCA}) or empirical models based on variational formulation; see \cite{isolaBone19,giorgio21,isolaBone12} among others.

\bibliographystyle{spbasic}
\bibliography{Ciallella}

\end{document}